\documentclass[12pt]{amsart}
\usepackage{amsmath,amssymb,amscd,verbatim}
\newcommand{\mc}{\mathcal}
\newcommand{\sub}{\subseteq}
\newcommand{\ol}{\overline}
\newcommand{\lra}{\Leftrightarrow}
\newcommand{\vo}{\varOmega}

\newcommand{\ra}{\Rightarrow}

\newcommand{\sm}{\setminus}
\DeclareMathOperator{\rad}{rad}

\newtheorem{theorem}{Theorem}[section]
\newtheorem{lemma}[theorem]{Lemma}
\newtheorem{proposition}[theorem]{Proposition}
\newtheorem{corollary}[theorem]{Corollary}
\newtheorem*{notation}{Notation}
\newtheorem{remark}[theorem]{Remark}

\theoremstyle{definition}
\newtheorem*{definition}{Definition}

\begin{document}
\title[Simple Overrings]{Integral Domains whose Simple Overrings are
Intersections of Localizations}

\author{Marco Fontana}
\address{Dipartimento di Matematica\\ Universit\`{a} degli Studi Roma Tre\\
    Largo San L. Murialdo, 1\\
    00146 Roma, Italy}
\email{fontana@mat.uniroma3.it}
\author{Evan Houston}
\address[Evan Houston and Thomas Lucas]{Department of Mathematics\\
    University of North Carolina at Charlotte\\
    Charlotte, NC 28223 U.S.A.}
\email[Evan Houston]{eghousto@email.uncc.edu}
\author{Thomas Lucas}
\email[Thomas Lucas]{tglucas@email.uncc.edu} 
\thanks{All three authors acknowledge support from the Cultural
Co-operation Agreement between Universit\`{a} degli Studi Roma Tre
and The University of North Carolina at Charlotte; the first-named
author was also partially supported by research grant MIUR
2001/2002 (Cofin 2000-MM01192794).}

\date{}
\begin{abstract} Call a domain $R$ an sQQR-domain if each
simple overring of $R$, i.e., each ring of the form $R[u]$ with
$u$ in the quotient field of $R$, is an intersection of
localizations of $R$. We characterize Pr\"ufer domains as
integrally closed sQQR-domains. In the presence of certain
finiteness conditions, we show that the sQQR-property is very
strong; for instance, a Mori sQQR-domain must be a Dedekind
domain.  We also show how to construct sQQR-domains which have
(non-simple) overrings which are not intersections of
localizations.
\end{abstract}
\maketitle


\section*{Introduction} \label{s:intro}

Throughout this work, $R$ will denote an integral domain with
quotient field $K$.  The Kaplansky transform of an ideal $I$ of
$R$ is denoted by $\vo(I)$ and is defined by $$\vo(I)=\{x \in K
\mid \text{for each } a \in I \text{ there is an integer } n(a)
\ge 1 \text{ such that }a^{n(a)}x \in R\}.$$ (The notation
$\vo_R(I)$ will be used when the context involves more than one
ring.)  In \cite{fh} the first- and second-named authors studied
domains $R$ each of whose overrings is a Kaplansky transform of an
ideal $I$ of $R$. This work is in part a sequel to that paper.  It
turns out, however, that our investigations depend more heavily on
the notions of unique minimal overrings and QQR-domains developed
by Gilmer and Heinzer in \cite{gh}.  Recall that a domain $R$ is a
\emph{QQR-domain} if each overring of $R$ is an intersection of
localizations of $R$. Davis \cite{da} showed that a Pr\"ufer
domain must have the QQR-property and asked whether the converse
is true. In their paper, Gilmer and Heinzer showed that the
converse does not hold, and they explored in depth the relation
between the Pr\"ufer and QQR-properties.

We review basic definitions and prove basic facts in
Section~\ref{s:basics}. In particular, we dub the domains of the
title \emph{sQQR-domains}, and we characterize Pr\"ufer domains as
being integrally closed sQQR-domains. The remainder of the paper
is devoted to studying the non-integrally closed case. Since
QQR-domains were characterized in \cite{gh}, we are interested in
sQQR-domains which are not QQR-domains.  Much of
Section~\ref{s:results} is concerned with showing that, in the
presence of many frequently studied finiteness conditions, there
are no such examples.  For instance, we show that a Mori
sQQR-domain is a Dedekind domain, that an sQQR-domain with
Noetherian spectrum is a Pr\"ufer domain, and that a semilocal
sQQR-domain with treed spectrum is a QQR-domain. In
Section~\ref{s:examples} we do provide many (non-local) examples
of sQQR-domains which are not QQR-domains. Finally, in the last
section, we briefly discuss connections with seminormality.


\section{Definitions and basic facts} \label{s:basics}

We begin by collecting some of the results we shall need. Recall
that for an ideal $I$ of $R$, the $v$-closure of $I$ is defined by
$I_v=(I^{-1})^{-1}$ and the $t$-closure by $I_t=\bigcup J_v$,
where the union is taken over all finitely generated subideals $J$
of $I$.  We assume familiarity with these star operations.

\begin{proposition} \label{p:basics} Let $I$ be an ideal of $R$.
Then:
\begin{itemize}
\item[(1)] $\vo(I)=\bigcap \{R_P \mid P \in \textup{Spec}(R), I
\nsubseteq P\}$.
\item[(2)] If $I$ is finitely generated and $S$ is a
multiplicatively closed subset of $R$, then $\vo_R(I)R_S =
\vo_{R_S}(IR_S)$.
\item[(3)] $\vo(I)=\vo(\rad I)$.
\item[(4)] If $I$ has the property that for each $P \in
\textup{Spec}(R)$ with $I \sub P$, we have $\vo(I) \nsubseteq
R_P$, then $\rad I$ is largest ideal $J$ for which
$\vo(I)=\vo(J)$.
\item[(5)] $\vo(I)=\vo(I_t)$.
\end{itemize}
\end{proposition}
\begin{proof} Statement (1) is proved in \cite{h}.  For statement
(2), one sees easily that the Kaplansky transform coincides with
the familiar Nagata transform for finitely generated ideals, and
the Nagata transform is known to localize well.  That
$\vo(I)=\vo(\rad I)$ is \cite[Lemma 3.1(c)]{f}.  For (4) suppose
that $I \sub P$ and that $\vo(J) =\vo(I) \nsubseteq R_P$.  By (1),
$J \sub P$.  Hence $J \sub \rad I$.  Finally, (5) appears in
\cite[Proposition 3.4]{f}.
\end{proof}

\begin{notation} \label{n:colon} Let $R$ be a domain with quotient field $K$.  If $U$
is a subset of $K$, we shall write $(R:U)$ for the fractional
ideal $\{z \in K \mid zU \sub R\}$ and $(R:_R U)$ for the ideal
$\{r \in R \mid rU \sub R\}$.
\end{notation}

\begin{lemma} \label{l:finiteext}
    Let $R$ be a domain, and let $U$ be a subset of $K$.  Then
\begin{itemize}
\item[(1)] $R[U] \subseteq \varOmega(R:_R U)$.
\item[(2)] If $U$ is finite, then $R[U]$ is the Kaplansky
transform of an ideal of $R$ $\lra$ $R[U]=\vo(R:_R U)$.
\end{itemize}
\end{lemma}
\begin{proof} Statement (1) follows easily from the definitions (with the convention
that $\vo(0)=K$).  Suppose that $U$ is finite and that
$R[U]=\vo(I)=\bigcap_{I \nsubseteq P} R_P$ for some ideal $I$ of
$R$. Let $P$ be a prime ideal of $R$ with $I \nsubseteq P$. Then
$U \subseteq R_{P}$, whence $(R:_R U) \nsubseteq P$.  Thus
$\varOmega(R:_R U) = \bigcap_{(R:_R U) \nsubseteq Q} R_{Q}
\subseteq R_{P}$. Hence $\varOmega(R:_R U) \subseteq \bigcap_{I
\nsubseteq P} R_{P} = R[U]$.  This proves (2).
\end{proof}

\begin{lemma} \label{l:same} Let $R$ be a domain with quotient field
$K$.  Then a finitely generated overring of $R$ is a Kaplansky
transform of an ideal of $R$ if and only if it is an intersection
of localizations of $R$.
\end{lemma}
\begin{proof} Let $U$ be a finite subset of $K$, and suppose that
there is a set of primes $\mc P$ in $R$ with $R[U]=\bigcap_{P \in
\mc P} R_P$.  For each $P \in \mc P$, we then have $(R:_R U)
\nsubseteq P$.  Hence $$R[U]=\bigcap_{P \in \mc P} R_P \supseteq
\bigcap_{(R:_R U) \nsubseteq Q} R_Q = \vo(R:_R U) \supseteq
R[U],$$ the last inclusion following from Lemma~\ref{l:finiteext}.
Thus $R[U] = \vo(R:_R U)$. The converse holds in general (without
the finiteness assumption) by Proposition~\ref{p:basics}(1).
\end{proof}

\begin{definition} \label{d:sqqr}  We say that a domain $R$ is an
\emph{sQQR-domain} (respectively, \emph{fQQR-domain}) if each
simple overring (respectively, each finitely generated overring)
of $R$ is an intersection of localizations of $R$.
\end{definition}

Recall from \cite{gh} that a domain is a \emph{QQR-domain} if each
overring of $R$ is an intersection of localizations of $R$.  In
\cite{fh} the authors defined an \emph{$\vo$-domain} to be one for
which each overring is a Kaplansky transform; hence a $\vo$-domain
is a QQR-domain. It was shown in \cite{fh} that the two properties
are not the same. However, if one restricts to finitely generated
overrings, then the properties are the same, as is shown by
Lemma~\ref{l:same}. Since \cite{gh} is a much older paper than
\cite{fh}, we have chosen to use the terminology ``sQQR''
(respectively, ``fQQR'') instead of ``s$\vo$'' (respectively,
``f$\vo$'').

In the integrally closed case, the sQQR- and fQQR-properties
coincide (Proposition~\ref{p:prufer}), but we have not been able
to determine whether these properties coincide in general. We can
easily see, however, that the sQQR- and $\vo$-properties are
distinct (even in the integrally closed case). This follows from
the fact that a valuation domain is automatically a QQR-domain
(and therefore an sQQR-domain), whereas a valuation domain need
not be an $\vo$-domain \cite[Example 2.16]{fh}.

In Section~\ref{s:examples}, we produce examples of sQQR-domains
which are not QQR domains.  However, we have not been able to
produce local examples of sQQR-domains which are not QQR. Indeed,
we devote much of the next section to showing that such examples
cannot exist in the presence of certain finiteness assumptions.
\bigskip


\section{Finiteness results} \label{s:results}

We have the following result in the integrally closed case.

\begin{proposition} \label{p:prufer}
    Let $R$ be integrally closed.  Then the following statements are equivalent.
\begin{itemize}
\item[(1)]  $R$ is an sQQR-domain.
\item[(2)]  $R$ is an fQQR-domain.
\item[(3)]  $R$ QQR-domain.
\item[(4)]  $R$ is a Pr\"{u}fer domain.
\end{itemize}
\end{proposition}
\begin{proof} Davis \cite[p. 197]{da} proved (4) $\ra$ (3).  Hence
it suffices to prove (1) $\ra$ (4).  This is a familiar argument
(cf. \cite[proof of Theorem 26.2]{g}).  Let $u \in K$, and
consider the simple overring $R[u^2]$. By hypothesis, this ring is
an intersection of localizations of $R$ and is therefore
integrally closed. Hence $u \in R[u^2]$, and so $u$ satisfies a
polynomial $f$ over $R$ such that some coefficient of $f$ is a
unit of $R$.  It now follows from the $u,u^{-1}$-lemma (as in
\cite[Lemma 19.14]{g}) that $R$ is a Pr\"{u}fer domain.
\end{proof}

\begin{proposition} \label{p:nonunits} Let $R$ be an sQQR-domain
which is not integrally closed, let $x \in \ol R \sm R$, and let
$\mc P$ denote the set of primes of $R$ which do not contain
$(R:_R x)$. Then
\begin{itemize}
\item[(1)] $\bigcup_{P \in \mc P} P$ is the set of nonunits of $R$, and
\item[(2)] $\mc P$ is infinite (so that $R$ has infinitely many
primes).
\end{itemize}
\end{proposition}
\begin{proof} If $s \in R$ but $s \notin \bigcup_{P \in \mc P} P$,
then $s$ is a unit of $\bigcap_{P \in \mc P} R_P=R[x]$ (the
equality following from Lemma~\ref{l:finiteext}). Since $x$ is
integral over $R$, this implies that $s$ is a unit of $R$. This
proves (1). For (2), note that (1) implies that $(R:_R x) \sub
\bigcup_{P \in \mc P} P$; prime avoidance then implies that $\mc
P$ must be infinite.
\end{proof}

This leads to our first restriction on (possible) local examples
of sQQR-domains which are not QQR-domains.

\begin{corollary} \label{c:locmaxnotprinc} Let $(R,M)$ be a local
sQQR-domain which is not integrally closed.  Then
\begin{itemize}
\item[(1)] $\dim R \ge 2$, and
\item[(2)] $M$ is not minimal over a principal ideal.
\end{itemize}
\end{corollary}
\begin{proof} Statement (1) is clear from
Proposition~\ref{p:nonunits}(2).  Statement (2) follows from
Proposition~\ref{p:nonunits}(1), upon observing that for a nonzero
element $a \in M$, $(R:_R 1/a)=(a)$.
\end{proof}

\begin{proposition} \label{p:loc} If $R_M$ is an sQQR-domain
(respectively, fQQR-domain) for each maximal ideal $M$ of $R$,
then $R$ is an sQQR-domain (respectively, fQQR-domain).
\end{proposition}
\begin{proof} For $x \in K$, we have \[R[x]=\bigcap_{M \in
\text{Max}(R)} (R[x])_{R \sm M} = \bigcap_{M \in \text{Max}(R)}
R_M[x].\]  For $M \in \text{Max}(R)$, since $R_M$ is an
sQQR-domain, $R_M[x]$ is an intersection of localizations of $R_M$
and therefore also of $R$. This takes care of the sQQR-property;
the proof for the fQQR-property is similar.
\end{proof}

Despite Proposition~\ref{p:loc}, neither the sQQR-property nor the
fQQR-property is a local property.  (We produce a specific example
of this in Section~\ref{s:examples}.) In fact, in
Proposition~\ref{p:locfindim} below, we show that a domain which
is locally a finite dimensional sQQR-domain is actually a Pr\"ufer
domain (and hence a QQR-domain).

We are able to localize the sQQR-property in some cases.  The
first such case is described in the following result; other cases
will be examined in Remark~\ref{r:locsqqr} below.

\begin{proposition} \label{p:treed} Let $R$ be a semilocal sQQR-domain
with $\textup{Spec}(R)$ treed.  Then $R_M$ is an sQQR-domain for
each maximal ideal $M$ of $R$.
\end{proposition}
\begin{proof} Fix a maximal ideal $M$.  Let $x \in K \sm M$, and let
$\mc P$ (respectively, $\mc Q$) denote the set of primes which do
not contain $(R:_R x)$ and which are contained in $M$
(respectively, are not contained in $M$).  Then we have
$R[x]=(\bigcap_{P \in \mc P} R_P) \bigcap (\bigcap_{Q \in \mc Q}
R_Q)$ by Proposition~\ref{p:basics}.  Set $S=R \sm M$.  Then
$$R_M[x]=R[x]_S = \left(\bigcap_{P \in \mc P} R_P\right)_S \bigcap
\left(\bigcap_{Q \in \mc Q} R_Q\right)_S=\left(\bigcap_{P \in \mc
P} R_P\right) \bigcap \left(\bigcap_{Q \in \mc Q} R_Q\right)_S.$$
We claim that $(\bigcap_{Q \in \mc Q} R_Q)_S \supseteq \bigcap_{P
\in \mc P} R_P$.  This will suffice to establish the result.  To
verify this, select a maximal ideal $N$ different from $M$.  Since
$\text{Spec}(R)$ is treed, Zorn's lemma produces a unique prime
ideal $N'$ which is maximal with respect to being contained in $M
\cap N$ (possibly, $N'=0$). If $(R:_R x) \nsubseteq N'$, then
$$\left(\bigcap_{Q \in \mc Q, Q \sub N} R_Q \right)_S \supseteq
(R_N)_S = R_{N'} \supseteq \bigcap_{P \in \mc P} R_P.$$  If $(R:_R
x) \sub N'$, then no $Q \in \mc Q$ satisfies $Q \sub N$, and the
intersection ``degenerates'' to $K$.  The claim now follows
easily.
\end{proof}

We require the concept of unique minimal overring: A proper
overring $T$ of a domain $R$ is said to be \emph{the unique
minimal overring} of $R$ if each proper overring of $R$ is
actually an overring of $T$. This concept was introduced in
\cite{gh} as a tool for studying QQR-domains.

\begin{proposition} \label{p:uniqueminovrg} Let $(R,M)$ be a local sQQR-domain, and suppose
that $R \subsetneqq \vo(M)$.  Then \begin{itemize}
\item[(1)] $\vo(M)$ is the unique minimal overring of $R$,
\item[(2)] $\vo(M)$ is an sQQR-domain, and
\item[(3)] if $\vo(M)=\ol R$, then $\ol R$ is a Pr\"ufer domain.
\end{itemize}
\end{proposition}
\begin{proof} (1)  Let $T$ be a proper overring of $R$, and let $t \in T \sm R$.  Then
$R[t]=\bigcap_{P \in \mc P} R_P$ for some set $\mc P$ of prime of
$R$ with $M \notin \mc P$.  Since $\vo(M)=\bigcap_{Q \ne M} R_Q$,
we have $\vo(M) \sub R[t] \sub T$.

(2) Pick $x \in K$, and write $R[x]=\bigcap_{P \in \mc P} R_P$ for
some set $\mc P$ of primes of $R$.  Since $\vo(M)$ is the unique
minimal overring of $R$, this yields $R[x]=\vo(M)[x]$; hence
$\vo(M)[x]$ is an intersection of localizations of $R$ and
therefore also of $\vo(M)$.

(3)  An integrally closed sQQR-domain is a Pr\"ufer domain by
Proposition~\ref{p:prufer}.
\end{proof}

\begin{proposition} \label{p:cases} Let $(R,M)$ be a local
sQQR-domain which is not integrally closed and for which $R
\subsetneqq (R:M)$. Then $(R:M)=(M:M)=\vo(M)$, $\vo(M)$ is the
unique minimal overring of $R$, and (exactly) one of the following
must occur:
\begin{itemize}
\item[(1)] $(R:M)$ is local with maximal ideal $M$,
\item[(2)] $(R:M)$ is local with maximal ideal different from $M$, or
\item[(3)] $(R:M)$ has two maximal ideals $N_1,N_2$ with $N_1 \cap N_2 = M$.
\end{itemize}
\end{proposition}
\begin{proof} By Corollary~\ref{c:locmaxnotprinc}, $M$ is not principal.
Since $M$ is maximal, we then have $(R:M)=(M:M)$, and this is a
proper overring of $R$.  Since $(R:M) \sub \vo(M)$, this implies
that $(R:M)=\vo(M)$ and that $\vo(M)$ is the unique minimal
overring of $R$ by Proposition~\ref{p:uniqueminovrg}. That this
situation results in just the three cases mentioned follows from
\cite[Section 2]{gh}.
\end{proof}

The next result shows that in cases (1) and (3) above, $R$ is a
QQR-domain.  We recall that a prime ideal of a domain $D$ is
\emph{unbranched} if $P$ is the only $P$-primary ideal of $D$ and
that in a Pr\"ufer domain this is equivalent to having $P$ be the
union of the prime ideals properly contained in $P$ \cite[Theorem
23.3]{g}.

\begin{proposition} \label{p:localsqqr2} Let $(R,M)$ be a local
sQQR-domain which is not integrally closed and for which $R
\subsetneqq (R:M)$. If $(R:M)$ is local with maximal ideal $M$ or
if $(R:M)$ has two maximal ideals, then
\begin{itemize}
\item[(1)] $(R:M)=\ol{R}$,
\item[(2)] $(R:M)$ is a Pr\"ufer domain, and
\item[(3)] $R$ is a QQR-domain.
\end{itemize}
\end{proposition}
\begin{proof}  We first deal with the case where $(R:M)$ has
maximal ideal $M$.  By \cite[Proposition 2.6]{gh}, $\ol R$ is a
valuation domain with maximal ideal $M$. Hence for $u \in \ol R$,
we have $uM \sub M$, and this implies that $(R:M)= \ol R$.  Thus
$\ol R = \vo(M) = \bigcap_{P \ne M} R_P = \bigcap_{Q \in \mc Q}
\ol R_Q$ for some set $\mc Q$ of primes of $\ol R$, with $Q \cap R
\ne M$ for each $Q \in \mc Q$. Since $\ol R$ is a valuation
domain, this implies that $M$ is unbranched in $\ol R$. By
\cite[Theorem 3.3]{gh}, $R$ is a QQR-domain.

Now suppose that $(R:M)$ has two maximal ideals $N_1$ and $N_2$.
In this case, $\ol R$ is a Pr\"ufer domain, $M$ is an ideal of
$\ol R$ and $N_1 \cap N_2=M$ by \cite[Proposition 2.5]{gh}. Hence
for $u \in \ol R$, we have $Mu \sub M \ol R \sub N_1 \cap N_2 =
M$,  and we have $u \in (R:M)$. Hence $\ol R=(R:M)$. We claim that
each $N_i$ is unbranched.  As before we have $\ol R = \bigcap_{P
\ne M} R_P = \bigcap_{Q \in \mc Q} \ol R_Q$ for some set $\mc Q$
of primes in $\ol R$, with $Q \cap R \ne M$ for each $Q$. Hence
for $Q \in \mc Q$, we must have $Q \ne N_i$ for $i=1,2$. If $N_i$
is branched, we can choose $t \in N_i$ but with $t$ in no other
prime of $\ol R$.  However, this yields $1/t \in \bigcap_{Q \in
\mc Q} \ol R_Q = \ol R$, a contradiction. Thus each $N_i$ is
unbranched, and $R$ is a QQR-domain, again by \cite[Theorem
3.3]{gh}.
\end{proof}

Recall that a \emph{pseudo-valuation domain} is a local domain
$(R,M)$ such that $(R:M)$ is a valuation domain with maximal ideal
$M$. By Proposition~\ref{p:localsqqr2}, we have the following:

\begin{corollary} \label{c:pvd}
    Let $(R,M)$ be a PVD.  Then $R$ is sQQR-domain $\Leftrightarrow$
$R$ is a QQR-domain.
\end{corollary}

We are not able to rule out the possibility of an example as in
case (2) of Proposition~\ref{p:cases}.  The following result
places some restrictions on such an example.

\begin{proposition} \label{p:case2} Let $(R,M)$ be a local
sQQR-domain which is not a QQR-domain.  Suppose that $R$ is as in
case (2) of Proposition~\ref{p:cases}. Set $T=(R:M)$, and let $N$
denote the maximal ideal of $T$.  Then $T$ is non-integrally
closed sQQR-domain, and $(T:N)=T=\vo_T(N)$.
\end{proposition}
\begin{proof}  $T$ is an sQQR-domain by
Proposition~\ref{p:uniqueminovrg}.  Suppose, by way of
contradiction, that $T$ is integrally closed.  Then $T$ is a
valuation domain by Proposition~\ref{p:prufer}.  In particular, we
must have $T=\ol R$. The fact that $T=\bigcap_{P \ne M} R_P$ then
implies that $T=\bigcap_{Q \in \mc Q} T_Q$ for some set $\mc Q$ of
primes of $T$ with $N \notin \mc Q$.  This, in turn, implies that
$N$ is unbranched. However, $M$ is an ideal of $T$ and is
therefore an $N$-primary ideal of $T$ different from $N$, a
contradiction. Thus $T$ is not integrally closed.  Now note that
$(T:M)=(R:M^2) \sub \vo(M) = T$, from which it follows immediately
that $(T:N)=T$. Moreover, $\vo(N)=\bigcap_{Q \nsupseteq N}T_Q
=\bigcap_{Q \nsupseteq MT} T_Q=\bigcap_{P \nsupseteq
M}R_P=\vo_R(M)=T$.
\end{proof}

\begin{corollary} \label{c:possiblelocalexample} If there is a local
sQQR domain $(R,M)$ which is not a QQR-domain, then there is such
an example with $(R:M)=R=\vo(M)$.
\end{corollary}

\begin{proposition} \label{p:locfindim} Let $R$ be a domain
such that $R_M$ is a finite dimensional sQQR-domain for each
maximal ideal $M$ of $R$. Then $R$ is a Pr\"ufer domain.
\end{proposition}
\begin{proof} If $R$ is a local sQQR-domain of dimension 1, then $R$ is a
valuation domain by Corollary~\ref{c:locmaxnotprinc} and
Proposition~\ref{p:prufer}. Let $(R,M)$ be a local sQQR-domain of
finite dimension greater than 1.  By induction, we may assume that
$R_P$ is a valuation domain for each nonmaximal ideal $P$. If
$R=\vo(M)=\bigcap_{P \ne M} R_P$, then $R$ is an integrally closed
sQQR-domain and is therefore a valuation domain by
Proposition~\ref{p:prufer}. Otherwise, $\vo(M)$ is the unique
minimal overring of $R$ by Proposition~\ref{p:uniqueminovrg}, and,
again since each $R_P$ is integrally closed, we have $\vo(M)=\ol
R$.  Hence $\ol R$ is a Pr\"ufer domain by
Proposition~\ref{p:uniqueminovrg}. Moreover, $\ol R$ has at most
two maximal ideals \cite[Corollary 2.2]{gh} and therefore $\ol R$,
and hence also $R$, has only finitely many prime ideals. It now
follows from Proposition~\ref{p:nonunits} that $R=\ol R$.  This
takes care of local $R$. The general case follows easily.
\end{proof}

We are now able to consider what happens when various finiteness
conditions are imposed on an sQQR-domain.  Recall that a domain
$R$ is said to be \emph{$v$-coherent} if for each finitely
generated ideal $I$ of $R$, we have $I^{-1}=J_v$ for some finitely
generated fractional ideal $J$ of $R$ and that $R$ is a
\emph{finite conductor domain} if each conductor ideal $(R:_R x)$,
where $x \in K$, is finitely generated.

\begin{proposition} \label{p:vcoh} If $R$ is an sQQR-domain which
is locally finite dimensional, and if, furthermore, $R$ is
$v$-coherent or a finite conductor domain, then $R$ is a Pr\"ufer
domain. In particular, a Noetherian sQQR-domain is a Dedekind
domain.
\end{proposition}
\begin{proof} Let $x \in K$.  If $R$ is a finite conductor domain,
then $(R:_R x)$ is finitely generated, and $R[x]$ is the Kaplansky
transform of a finitely generated ideal. If $R$ is $v$-coherent,
then $(R:_R x)=(R:(1,x))=I_t$ for some finitely generated ideal
$I$. Hence by Proposition~\ref{p:basics} and
Lemma~\ref{l:finiteext}, we again have that $R[x]$ is the
Kaplansky transform of a finitely generated ideal. It then follows
from Proposition~\ref{p:basics} that each localization of $R$ is
also an sQQR-domain. The result now follows from
Proposition~\ref{p:locfindim}.
\end{proof}

Recall that a domain $R$ is said to have \emph{Noetherian
spectrum} if it satisfies the ascending chain condition on radical
ideals and to be a \emph{Mori domain} if it satisfies the
ascending chain condition on divisorial ideals.

\begin{remark} \label{r:locsqqr} From Proposition~\ref{p:loc} and the proof of
Proposition~\ref{p:vcoh}, we obtain that if $R$ is a $v$-coherent
or a finite conductor domain, then $R$ is sQQR if and only if $R$
is locally sQQR.  We shall see later (in the proof of
Theorem~\ref{t:noespecsqqr}) that the sQQR-property is also a
local property for domains with Noetherian spectrum and for Mori
domains.
\end{remark}

\begin{lemma} \label{l:accraddiv} Let $(R,M)$ be a non-integrally
closed local sQQR-domain which satisfies the ascending chain
condition on radicals of conductor ideals.  Then $\vo(M)$ is the
unique minimal overring of $R$.
\end{lemma}
\begin{proof} Since $R$ is an sQQR-domain, for each $x \in K$ we
have $R[x]=\vo(R:_R x)=\vo(\rad(R:_R x))$ by
Proposition~\ref{p:basics} and Lemma~\ref{l:finiteext}. Moreover,
the ideal $(R:_R x)$ satisfies the hypothesis of
Proposition~\ref{p:basics}(4); hence for $x,y \in K$ we have $R[x]
\subsetneqq R[y]$ if and only if $\rad(R:_R y) \subsetneqq
\rad(R:_R x)$.  It then follows by our assumption that we may pick
$u \in \ol R \sm R$ such that $R[u]$ is a minimal extension of $R$
(i.e., there are no rings properly between $R$ and $R[u]$). By
\cite[Lemma 2.3]{gh}, $M$ is the conductor of $R$ in $R[u]$,
whence $R \subsetneqq (R:M) \sub \vo(M)$.  The lemma now follows
from Proposition~\ref{p:uniqueminovrg}.
\end{proof}

It is known that in a Mori domain the radical of a divisorial
ideal is again divisorial. Hence domains with Noetherian spectrum
and Mori domains satisfy the ascending chain condition on radicals
of conductor ideals.

\begin{theorem} \label{t:morisqqr} If $R$ is a Mori sQQR-domain, then $R$
is a Dedekind domain.
\end{theorem}
\begin{proof}
We first handle the case where $R$ is actually assumed to be a
Mori QQR-domain.  In this case, since both properties localize, we
may as well assume that $R$ is also local with maximal ideal $M$.
If $R$ is integrally closed, then $R$ is a valuation domain, and
it is well known that a Mori Pr\"ufer domain is a Dedekind domain
(see, e.g., the proof of \cite[Proposition 2.6]{cgh}). Hence we
also assume that $R$ is not integrally closed.  By \cite[Theorem
3.3]{gh}, $\ol R$ is a Pr\"ufer domain, and the (at most two)
maximal ideals of $\ol R$ are unbranched.  However, if $Q$ is a
nonzero, nonmaximal prime ideal of $\ol R$, then by \cite[Lemma
3.4]{fh}, $R_{Q \cap R}=\ol R_Q$, whence $R_{Q \cap R}$ is a Mori
valuation domain.  It is well known that this implies that $\ol
R_Q$ is a rank one discrete valuation domain, whence $Q$ has
height one. Thus the maximal ideals of $\ol R$ have finite height,
which contradicts the fact that they are unbranched.

For the rest of the proof, we explicitly assume that $R$ is not a
QQR-domain.  Since for an ideal $I$ of a Mori domain, we have
$I_t=I_v=J_v=J_t$ for some finitely generated ideal $J$ of $R$, we
have $\vo(I)=\vo(J)$ by Proposition~\ref{p:basics}, so that the
sQQR-property localizes. Hence we may assume that $R$ is a local
Mori sQQR-domain with maximal ideal $M$. Of course, we may also
assume that $R$ is not integrally closed.  By
Lemma~\ref{l:accraddiv}, $\vo(M)=(R:M)$ is the unique minimal
overring of $R$.  Hence $\vo(M)$ is integral over $R$, and by
\cite[Lemma 2.3]{gh}, $M$ is the conductor of $R$ in $\vo(M)$.  In
particular, $M$ is divisorial.  By Propositions~\ref{p:cases},
\ref{p:localsqqr2}, and \ref{p:case2}, $T=\vo(M)$ is a local
sQQR-domain which is not integrally closed, and the maximal ideal
$N$ of $T$ satisfies $(T:N)=T$, so that $N$ is nondivisorial.
Moreover, since $\vo(M)=(R:M)$, $\vo(M)$ is again a Mori domain by
\cite[I, Th\'eor\`eme 2]{r}. However, as we just showed, this
implies that the maximal ideal $N$ of $T$ is divisorial, a
contradiction.
\end{proof}

\begin{theorem} \label{t:noespecsqqr}  Let $R$ be an sQQR-domain
with Noetherian spectrum.  Then $R$ is a Pr\"ufer domain.
\end{theorem}
\begin{proof} A ring with Noetherian spectrum has the property that
each radical ideal is the radical of a finitely generated ideal.
Hence each Kaplansky transform is the transform of a finitely
generated ideal.  It then follows from Proposition~\ref{p:basics}
that the sQQR-property localizes for domains with Noetherian
spectrum. Hence we may assume that $R$ is a local sQQR-domain with
maximal ideal $M$.  We may also assume that $R$ is not integrally
closed. If $R$ is a QQR-domain, then $M$ is unbranched by
\cite[Theorem 3.3]{gh}. However, $M$ is the radical of a finitely
generated ideal, which contradicts that fact that $M$ is
unbranched.  (To see this, let $M = \rad J$, $J$ finitely
generated.  Since $M$ is unbranched, we must have $M=J$, i.e., $M$
is finitely generated.  However, Nakayama's lemma then guarantees
that $M \ne M^2$, and $M^2$ is then an $M$-primary ideal distinct
from $M$, a contradiction.)  Hence we may assume that $R$ is not a
QQR-domain. By Lemma~\ref{l:accraddiv}, $\vo(M)$ is the unique
minimal overring of $R$, and, as in the proof of
Theorem~\ref{t:morisqqr}, $M$ is the conductor of $R$ in $\vo(M)$,
so that $M$ is divisorial. Also as in the proof of
Theorem~\ref{t:morisqqr}, we obtain that $T=\vo(M)$ is a local
sQQR-domain with nondivisorial maximal ideal $N$.  However, it is
easy to see that the contraction map from $\text{Spec}(T)$ to
$\text{Spec}(R)$ is a homeomorphism, from which it follows that
$T$ also has Noetherian spectrum.  But then the argument just
given shows that $N$ is divisorial, a contradiction.
\end{proof}

A one dimensional domain in which each nonzero element is
contained in only finitely many primes has Noetherian spectrum.
Hence we have the following result.

\begin{corollary} \label{c:onedimprufer}
    Let $R$ be a one-dimensional domain, and assume that each nonzero element of $R$
is contained in only finitely many primes. Then $R$ is an
sQQR-domain $\lra$ $R$ is a Pr\"ufer domain.
\end{corollary}

\begin{theorem} \label{t:semilocsqqr} If $R$ is a semilocal
sQQR-domain with treed spectrum, then $R$ is a QQR-domain.
\end{theorem}
\begin{proof}  By \cite[Theorem 1.9]{gh}, it suffices to prove
that each localization of $R$ is a QQR-domain.  Hence by
Proposition~\ref{p:treed}, we may assume that $R$ is local with
maximal ideal $M$.  We may also assume that $R$ is not integrally
closed.  Let $x \in \ol R \sm R$, and write $R[x]=\bigcap_{(R:_R
x) \nsubseteq P} R_P$.  Since $x$ is integral over $R$, this
intersection cannot contain $R_Q$ for any prime $Q \ne M$.  Thus
we must have $\rad(R:_R x) = M$.  Hence $R[x] = \vo(M)$.  By
Proposition~\ref{p:uniqueminovrg}, $\vo(M)$ is the unique minimal
overring of $R$, and, since $x$ was chosen arbitrarily, we have
$\vo(M)=\ol R$, so that $\ol R$ is a Pr\"ufer domain.  Note that
in this situation $\ol R=\bigcap_{P \ne M}R_P$, and so $M$ must be
the union of the chain of primes properly contained in $M$.  We
claim that each maximal ideal of $\ol R$ is unbranched. If not,
let $N$ be a branched maximal ideal of $\ol R$.  Then, since $\ol
R$ is a Pr\"ufer domain, $N=\rad \ol Ru$ for some $u \in \ol R$
\cite[Theorem 23.3]{g}.  By \cite[Section 2]{gh}, $u^2 \in R$,
whence $M = \rad Ru^2$. However, this contradicts the fact that
$M$ is the union of the chain of primes properly contained in $M$.
Hence $R$ is a QQR-domain by \cite[Theorem 3.3]{gh}.
\end{proof}


\section{Examples} \label{s:examples}

In this section, we construct examples of sQQR-domains which are
not QQR-domains.

\begin{lemma} \label{l:omegapullback}
    Consider the following pullback diagram.
\[
    \begin{CD}
        R   @>>>    D   \\
        @VVV        @VVV    \qquad \qquad   \\
        T   @>\varphi>>    k   \\
    \end{CD}
\]
Here, $T$ is a domain, $k=T/M$ for some maximal ideal $M$ of $T$,
$\varphi$ is the canonical homomorphism, and $D$ is a domain
contained in $k$. Let $I$ be an ideal of $R$ which properly
contains $M$. Then $\varOmega_{R}(I) =
\varphi^{-1}(\varOmega_{D}(\varphi(I)))$.
\end{lemma}
\begin{proof}
    Let $x \in \varOmega_{R}(I)$.  Note that since $I \supsetneqq M$
we have $x \in T$.  Let $u=\varphi(a) \in \varphi(I)$, where $a
\in I$.  Then $xa^{n} \in R$ for some positive integer $n$.  It
follows that $\varphi(x)\varphi(a)^{n} \in D$.  Thus $x \in
\varphi^{-1}(\varOmega_{D}(\varphi(I)))$.  Now let $y \in
\varphi^{-1}(\varOmega_{D}(\varphi(I)))$, and let $b \in I$.  We
have $\varphi(y) \in \varOmega_{D}(\varphi(I))$, so that
$\varphi(yb^{m}) = \varphi(y)\varphi(b)^{m} \in D$ for some
positive integer $m$; that is, $yb^{m} \in R$.  It follows that $y
\in \vo_R(I)$.
\end{proof}

\begin{proposition} \label{p:sqqrpullback}
    Consider the pullback diagram of Lemma~\ref{l:omegapullback},
and assume that $T=V$ is a valuation domain and that $k$ is the
quotient field of $D$.  Then $R$ is an sQQR-domain (respectively,
fQQR-domain) $\lra$ $D$ is an sQQR-domain (respectively,
fQQR-domain).
\end{proposition}
\begin{proof}
    We give the details for the sQQR-case, the fQQR-case being
similar.  Assume that $R$ is an sQQR-domain. Let $u \in k$, and
let $I = \varphi^{-1}(D:_D u)$.  Pick $v \in V$ with
$\varphi(v)=u$. Then  $I = (R:_R v)$, and $R[v] =
\varOmega_{R}(I)$.  Hence by Lemma~\ref{l:omegapullback}, $D[u] =
\varphi(R[v]) = \varphi(\varOmega_{R}(I)) =
\varOmega_{D}(\varphi(I))) = \varOmega_{D}(D:_D u)$.  Thus $D$ is
an sQQR-domain.

For the converse, let $x \in K$.  If $x \in V$ we can use the same
techniques to show that $R[x] = \varOmega_{R}(R:_R x)$.  If $x
\notin V$, then $x^{-1} \in M \subseteq R$, and it is well known
(and easy to show) that $R[x] = \varOmega_{R}(x^{-1}R)$.
\end{proof}

\begin{proposition} \label{p:notqqr} Let $T$ be a Pr\"ufer domain with the following
properties: \begin{itemize}
\item[(1)] $T$ is one dimensional with $\textup{Max}(T)=\{M\} \cup \mc
P$, with $\mc P$ infinite,
\item[(2)] Each element of $M$ is in $P$ for almost all $P \in \mc P$ (equivalently,
each finitely generated ideal contained in $M$ is contained in $P$
for almost all $P \in \mc P$),
\item[(3)] Each element $a \in T \sm M$ satisfies $a \notin P$ for infinitely many
$P \in \mc P$, and
\item[(4)] $T/M$ admits a proper subfield $F$ such that there are no fields between
$F$ and $T/M$.
\end{itemize}

Then $R$ is an sQQR-domain but is not a QQR-domain, where $R$ is
defined by the following pullback diagram: \[
    \begin{CD}
        R   @>>>    F   \\
        @VVV        @VVV    \qquad \qquad   \\
        T   @>\varphi>>    T/M.   \\
    \end{CD}
\]
Moreover, $R_M$ is not an sQQR-domain.
\end{proposition}
\begin{proof} Let $K$ denote the common quotient field of $R$
and $T$.  Suppose that $\mc Q$ is an infinite subset of $\mc P$.
Then for $u \in \bigcap_{Q \in \mc Q} T_Q$, we have that the
finitely generated ideal $(T:_T u) \nsubseteq Q$ for each $Q \in
\mc Q$, whence by (2) $(T:_T u) \nsubseteq M$.  Thus $u \in T_M$.
That is, $\bigcap_{Q \in \mc Q} T_Q \sub T_M$ for each infinite
subset $\mc Q$ of $\mc P$.

Now let $x \in K \sm R$.  We consider two cases.

Case 1.  Suppose that $T[x] \sub T_M$.  Then $(T:_T x) \nsubseteq
M$, whence by (3) $(T:_T x) \nsubseteq P$ for infinitely many $P
\in \mc P$. Thus $T[x]=\bigcap_{Q \in \mc Q} T_Q$ for some
infinite subset $\mc Q$ of $\mc P$.  Since $(T:_T x) \nsubseteq
M$, we can choose $t \in T \sm M$ with $tx \in T$. Now $t \notin
M$ implies that $\varphi(t) \ne 0$, and we may find $t' \in T$
with $\varphi(tt')=1$.  Let $r=tt'$.  Then $r \in R \sm M$, and
$rx \in T$.  If $rx \notin R$, then, since $rx \in T$ and $R \sub
T$ is a minimal extension (there are no rings properly between $R$
and $T$), we have $T=R[rx] \sub R[x]$, whence
$R[x]=T[x]=\bigcap_{Q \in \mc Q} T_Q = \bigcap_{Q \in \mc Q} R_{Q
\cap R}$.  Thus we may assume that $rx \in R$, so that $(R:_R x)
\nsubseteq M$, and we have $R[x] \sub R_M$.  We shall show, in
fact, that $R[x] = T[x] \cap R_M$, which will complete the proof
in this case.  Note that $M[x]$ is a maximal ideal of $T[x]$, and
it follows that $M[x]$ is also a maximal ideal of $R[x]$.  Thus we
have the following pullback diagram: \[
    \begin{CD}
        R[x]   @>>>    R[x]/M[x]   \\
        @VVV        @VVV    \qquad \qquad   \\
        T[x]   @>>>    T[x]/M[x].
    \end{CD}
\]  Moreover, $R/M \sub R[x]/M[x] \sub R_M/MR_M$, whence $R[x]/M[x] = F$.  Similarly,
$T[x]/M[x] = T/M$.  Hence $R[x] \sub T[x]$ is a minimal extension.
However, we have $R[x] \sub T[x] \cap R_M \subsetneqq T[x]$,
whence $R[x] = T[x] \cap R_M$, as desired.

Case 2.  Suppose that $T[x] \nsubseteq T_M$.  Then $(T:_T x) \sub
M$, so that $(T:_T x) \sub P$ for almost all $P \in \mc P$ by (2).
Hence $T[x] = T_{P_1} \cap \cdots \cap T_{P_n}$ for some finite
subset $\{P_1, \ldots, P_n\}$ of $\mc P$. It follows that
$MT[x]=T[x]$. However, $MT[x]=M[x]=MR[x]\sub R[x]$, so that
$R[x]=T[x] = \bigcap_{j=1}^n T_{P_j} = \bigcap_{j=1}^n R_{P_j \cap
R}$. This completes the proof that $R$ has the sQQR-property.

Now by \cite[Theorem 1.9]{gh} the QQR-property is a local
property.  Moreover, $M$ is branched ($R$ is one dimensional),
whence by \cite[Theorem 3.3]{gh}, $R_M$ does not have the
QQR-property.  Hence $R$ is not a QQR-domain.  Finally, observe
that, since $R_M$ is not integrally closed,
Corollary~\ref{c:locmaxnotprinc} assures that $R_M$ is not an
sQQR-domain.
\end{proof}

There remains the construction of domains with the properties
described in Proposition~\ref{p:notqqr}.  Begin with an almost
Dedekind domain $T_0$ which is not Dedekind (so that there are
necessarily infinitely many maximal ideals), which has precisely
one non-finitely generated maximal ideal $M$, and which satisfies
property (2) of the proposition; \cite[Example 2, pp. 338-339]{g3}
is one such example. We claim that condition (3) of the
proposition is then automatically satisfied.  To see this pick $a
\in T_0 \sm M$; then $(M,a)=T_0$, and we may choose $c \in M$ so
that $(a,c)=T_0$. Since $c$ lies in almost all $P_i$, $a$ fails to
be in almost all $P_i$. As for condition (4), if it is not already
satisfied, replace $T_0$ by $T=T_0(X)$ ($T_0(X)=T_0[X]_S$, where
$S$ is the multiplicatively closed subset of $T_0[X]$ consisting
of all polynomials having unit content). By \cite[Proposition
36.7]{g}, $T$ is again an almost Dedekind domain.  Moreover, the
maximal ideals of $T$ are just the extensions of the maximal
ideals of $T_0$; hence $T$ has properties (1)-(3) above.  The
pertinent residue field is $T/MT(X)=(T_0/M)(X)$, which admits
appropriate subfields $F$, e.g., $F=(T_0/M)(X^2)$.

Since valuation domains are (s)QQR-domains,
Propositions~\ref{p:sqqrpullback} and \ref{p:notqqr} can be used
to produce examples of arbitrary dimension of sQQR-domains which
are not QQR-domains.


\section{Seminormality} \label{s:seminormality}

Recall that if $R$ is a domain with quotient field $K$, then $R$
is \emph{seminormal} if whenever $x \in K$ with $x^2,x^3 \in R$,
we have $x \in R$.  We show below that a QQR-domain is seminormal,
but we have not been able to determine whether an sQQR-domain (or
even an fQQR-domain) need be seminormal.  The example of an
sQQR-domain which is not a QQR-domain discussed at the end of
Section~\ref{s:examples} is seminormal (it is also an
fQQR-domain), so a seminormal sQQR-domain need not be a
QQR-domain, but we do not know whether a local seminormal
sQQR-domain (or fQQR-domain) must be a QQR-domain.

\begin{proposition} \label{p:qqrseminormal}  A QQR-domain is
seminormal.
\end{proposition}
\begin{proof} It suffices to establish the result locally.
Hence we assume that $(R,M)$ is a local QQR-domain, and we may as
well assume that $R$ is not Pr\"ufer.  By \cite[Section 2 and
Theorem 3.3]{gh} either $\ol R$ is a valuation domain with maximal
ideal $M$ or $\ol R$ is a Pr\"ufer domain with two maximal ideals
$N_1,N_2$ such that $N_1 \cap N_2 = M$.  In particular, $M$ is a
radical ideal of $\ol R$.  Now let $x \in K$ be such that $x^2,x^3
\in R$. Then $x \in \ol R$.  If $x^2 \notin M$, then $x=x^3/x^2
\in R$. If $x^2 \in M$, then since $M$ is a radical ideal of $\ol
R$, we have $x \in M \sub R$, as desired.
\end{proof}

We close with a characterization of seminormal fQQR-domains.

\begin{proposition} \label{p:fqqrseminormal} The following
statements are equivalent for an fQQR-domain $R$.
\begin{itemize}
\item[(1)] Each overring of $R$ is seminormal.
\item[(2)] Each integral finitely generated overring of $R$ is seminormal.
\item[(3)] $R[x]=R[x^2,x^3]$ for each $x \in K$.
\item[(4)] $\rad(R:_R x) = \rad(R:_R x^2) \cap \rad(R:_R x^3)$ for
each $x \in K$.
\item[(5)] $R$ is seminormal.
\end{itemize}
If $R$ satisfies any of these equivalent conditions, then $\ol R$
is a Pr\"ufer domain.
\end{proposition}
\begin{proof} The implications (3) $\lra$ (1) $\ra$ (2)
$\ra$ (5) are straightforward. Let $x \in K$. Then, using
Proposition~\ref{p:basics}, we have $R[x^2,x^3]=R[x]$ $\lra$
$\vo(R:_R (x^2,x^3)) = \vo(R:_R x)$.  As in the proof of
Lemma~\ref{l:accraddiv}, this last equality holds $\lra$
$\rad(R:_R (x^2,x^3) = \rad(R:_R x)$. Since $(R:_R (x^2,x^3)) =
(R:_R x^2) \cap (R:_R x^3)$, this proves (3) $\lra$ (4).  Assume
(5), and let $a \in (R:_R (x^2,x^3))$. Then $ax^2,ax^3 \in R$,
whence $(ax)^2,(ax)^3 \in R$. Since $R$ is seminormal, $ax \in R$.
Thus $\rad(R:_R x^2) \cap \rad(R:_R x^3) \sub \rad(R:_R x)$. The
converse always holds. Hence (5) $\ra$ (4).  Finally, the last
statement holds by \cite[Theorem 2.3]{adh}.
\end{proof}

\end{document}